\documentclass[12pt,  reqno]{amsart}

\newtheorem{theorem}{Theorem}
\theoremstyle{plain}

\newtheorem{corollary}[theorem]{Corollary}

\begin{document}
\title[] {On the $q$-Euler numbers and polynomials with weight $0$}

\author[] {T. Kim}
\address{Taekyun Kim. Division of General Education-Mathematics \\
Kwangwoon University, Seoul 139-701, Republic of Korea  \\}
\email{tkkim@kw.ac.kr}

\maketitle

{\footnotesize {\bf Abstract} \hspace{1mm} {The purpose of this
paper is to investigate some properties of $q$-Euler numbers and
polynomials with weight $0$. From those $q$-Euler numbers with weight $0$, we derive some identies on the $q$-Euler numbers and
polynomials with weight $0$. }

\section{Introduction}
Let $p$ be a fixed odd prime number. Throughout this paper
 $\Bbb Z_p$, $\Bbb Q_p$, and $\Bbb C_p$ will denote
the ring of $p$-adic rational integers, the field of $p$-adic rational
numbers, and the completion of algebraic closure of $\Bbb Q_p$.
The $p$-adic absolutely value is defined by
$|x|_p=1/p^{r}$ where  $x=p^{r} s/t$ with $(p,s)=(p,t)=(s,t)=1$ and $r \in \Bbb Q$.
In this paper, we assume that $\alpha \in \Bbb Q$ and $q \in \Bbb C_p$ with $|1-q|_p < 1$.
As well known definition, the Euler polynomials are defined by
$$ \frac {2}{e^t +1} e^{xt} = e^{E(x)t} = \sum_{n=0}^{\infty} E_{n} (x) \frac {t^n}{n!},$$
with the usual convention about
replacing ${E}^{n} (x)$ by ${E}_{n}(x)$ (see [1-15]).

In this special case, $x=0,  E_n (0)=E_n$ are called the $n$-th
 Euler numbers (see [1]).
Recently, the $q$-Euler numbers with weight $\alpha$ are defined by
\begin{eqnarray}
{{\tilde E}}_{0,q}^{(\alpha)} =1, \quad \text{and}
 \quad q( q^{\alpha}{\tilde E}_q^{(\alpha)} +1)^n +{\tilde E}_{n,q}^{(\alpha)}= 0\
\ \text{if}\ \ n>0, \end{eqnarray}
with the usual convention about
replacing ${( \tilde E_{q}^{(\alpha)} )}^n$ by $\tilde{E}_{n,q}^{(\alpha)}$ (see [3,12]).
The $q$-number of $x$ is defined by $[x]_q =\frac{1- q^x}{1-q}$ (see [1-15]).
Note that $\lim_{q \rightarrow 1} [x]_q =x$.
Let us define the notation of $q$-Euler numbers with weight $0$ as ${\tilde E_{n,q}}^{(0)} = {\tilde E_{n,q}}$.
The purpose of this paper is to investigate some interesting identities on the $q$-Euler numbers with weight $0$.

\medskip

\section{On the extended $q$-Euler numbers of higher-order with weight 0}

Let $C(\Bbb Z_p )$ be the space of continuous functions on $\Bbb Z_p$. For $f \in C(\Bbb Z_p )$, the fermionic $p$-adic
$q$-integral on $\Bbb Z_p$ is defined by Kim as follows :

\begin{eqnarray} I_{q}(f)&=&\int_{\Bbb Z_p
} f(x) d \mu_{-q} (x)\notag
\\ &=& \lim_{N \rightarrow \infty} \frac {[2]_q}{1+q^{p^N}}
\sum_{x=0}^{p^N -1} f(x)(-q)^x, \text{ (see [1-12])}.
\end{eqnarray}

By (2), we get
\begin{eqnarray}
q^n I_{q}(f_n)+(-1)^{n-1}I_{q}(f)=
[2]_q  \sum_{l=0}^{n -1}(-1)^{n-1-l}f(l)q^l,
\\
\text{where} \,\ f_{n} (x)=f(x+n)\,\ \quad \text{and} \quad n \in \Bbb N \,\  \text{ (see [4, 5])}.
 \notag
\end{eqnarray}

By (1), (2) and (3), we see that
\begin{eqnarray}
\int_{\Bbb Z_p} [x]_{q^\alpha}^n d\mu_{-q}(x) =
\tilde E_{n,q}^{(\alpha)}= \frac {[2]_{q}}{(1-q)^n [\alpha]_{q}^{n}}
 \sum_{l=0}^{n} \binom {n}{l} (-1)^{l} \frac {1}{1+q^{\alpha l+1}}
.
\end{eqnarray}

In the special case, $n=1$, we get
\begin{eqnarray}
\int_{\Bbb Z_p} e^{xt} d\mu_{-q}(x)=
\frac {[2]_{q}}{ qe^{t} +1 } = \frac {1+q^{-1}}{e^t +q^{-1}}=
 \sum_{n=0}^{\infty} H_n {(-q^{-1})} \frac {t^n}{n!},
\end{eqnarray}
 where $H_n (-q^{-1})$ are the $n$-th Frobenius-Euler numbers.
From (5), we note that the
$q$-Euler numbers with weight $0$ are given by
\begin{eqnarray}
\tilde E_{n,q} = \int_{\Bbb Z_p} x^{n} d\mu_{-q}(x) = H_n (-q^{-1}), \,\ \text{for} \quad n \in \Bbb Z_+
.
\end{eqnarray}

Therefore, by (6), we obtain the following theorem.
\begin{theorem}
For $n \in \Bbb Z_+$, we have
$$
\tilde E_{n,q} = H_n (-q^{-1}),
$$
where $H_n (-q^{-1})$ are called the $n$-th Frobenius-Euler numbers.
\end{theorem}

Let us define the generating function of the $q$-Euler numbers with weight $0$ as follows:
\begin{eqnarray}
\tilde F_q (t)&=& \sum_{n=0}^{\infty} \tilde E_{n,q}
\frac  {t^n}{n!}.
\end{eqnarray}

Then, by (4) and (7), we get
\begin{eqnarray}
\tilde F_q (t) = [2]_q \sum_{m=0}^{\infty} (-1)^m q^m e^{mt} =
\frac  {1+q}{qe^{t}+1}.
\end{eqnarray}

Now we define the $q$-Euler polynomials with weight $0$ as follows:
\begin{eqnarray}
\sum_{n=0}^{\infty} \tilde E_{n,q} (x)
\frac  {t^n}{n!}= \frac {1+q}{qe^{t}+1}e^{xt}.
\end{eqnarray}

Thus, (5) and (9), we get
\begin{eqnarray}
\int_{\Bbb Z_p} e^{(x+y)t} d\mu_{-q}(y)
= \frac {1+q}{qe^{t}+1}e^{xt}
=\sum_{n=0}^{\infty} \tilde E_{n,q} (x)
\frac  {t^n}{n!}.
\end{eqnarray}

From (10), we have
\begin{eqnarray}
\sum_{n=0}^{\infty} \tilde E_{n,q} (x) \frac  {t^n}{n!} = \left( \frac {1+q^{-1}}{e^{t}+q^{-1}} \right ) e^{xt}
=\sum_{n=0}^{\infty} H_n  (-q^{-1},x) \frac  {t^n}{n!}
,
\end{eqnarray}
where $H_n (-q^{-1},x)$ are called the $n$-th Frobenius-Euler polynomials (see[9]).

Therefore, by (11), we obtain the following theorem.
\begin{theorem}
For $n \in \Bbb Z_+$, we have
$$
\tilde E_{n,q}(x) = \int_{\Bbb Z_p} {(x+y)}^n  d\mu_{-q}(x)
= H_n (-q^{-1},x),
$$
where $H_n (-q^{-1},x)$ are called the $n$-th Frobenius-Euler polynomials.
\end{theorem}

From (3) and Theorem 2, we note that
\begin{eqnarray}
q^n H_m (-q^{-1},n)+ H_m
(-q^{-1}) = [2]_q \sum_{l=0}^{n-1} (-1)^l l^m q^l ,
\end{eqnarray}
where $n \in \Bbb N$  with  $n \equiv 1 \pmod 2$.

Therefore, by (12), we obtain the following corollary.
\begin{corollary}
For $n \in \Bbb N$, with  $n \equiv 1 \pmod 2$ and $m \in \Bbb Z_+$, we have
$$
q^n H_m (-q^{-1},n)+ H_m(-q^{-1}) = [2]_q \sum_{l=0}^{n-1} (-1)^l l^m q^l.
$$
\end{corollary}

In particular, $q=1$, we get $E_m (n)+E_m = 2 \sum_{l=0}^{n-1} (-1)^l l^{m}$,
where $E_m$ and $E_m (n)$ are called the $m$-th Euler numbers and polynomials which are defined by
$$
\frac {2}{e^t +1} = \sum_{m=0}^{\infty} E_m \frac {t^m}{m!} \quad  \text{and} \quad
\frac {2}{e^t +1}e^{xt} = \sum_{m=0}^{\infty} E_m (x) \frac {t^m}{m!}.
$$

By (3), we easily see that
\begin{eqnarray}
q \int_{\Bbb Z_p} f(x+1) d\mu_{-q}(x)+ \int_{\Bbb Z_p} f(x) d\mu_{-q}(x)
= [2]_{q}f(0).
\end{eqnarray}

Thus, by (13), we get
\begin{eqnarray*}
[2]_q
&=& q \int_{\Bbb Z_p} e^{(x+1)t} d\mu_{-q}(x)+ \int_{\Bbb Z_p} e^{xt} d\mu_{-q}(x) \\
&=& \sum_{n=0}^{\infty} \left ( q \int_{\Bbb Z_p} (x+1)^n  d\mu_{-q}(x)+ \int_{\Bbb Z_p} x^n  d\mu_{-q}(x) \right )
\frac {t^n}{n!} \\
&=&  \sum_{n=0}^{\infty} \left( q H_n (-q^{-1},1) + H_n (-q^{-1})\right ) \frac {t^n}{n!}.
\end{eqnarray*}

Therefore, by(13), we obtain the following theorem.
\begin{theorem}
For $n \in \Bbb Z_+$, we have
\begin{eqnarray*}
 q H_n (-q^{-1},1) + H_n (-q^{-1})
 = \left \{
\begin{array}{ll}
 1+q , \  \ & \hbox{if} \ \ n=0,
\vspace{2mm} \\ 0, \ \ & \hbox{if} \ \ n>0,
\end{array}\right.
\end{eqnarray*}
\end{theorem}

where $H_n (-q^{-1},x)$ are called the $n$-th Frobenius-Euler polynomials and
$H_n (-q^{-1})$ are called the $n$-th Frobenius-Euler numbers.
In particular, $q=1$, we have

\begin{eqnarray*}
 E_n (1)+E_n  = \left \{
\begin{array}{ll}
 2 , \  \ & \hbox{if} \ \ n=0,
\vspace{2mm} \\ 0, \ \ & \hbox{if} \ \ n>0,
\end{array}\right.
\end{eqnarray*}
where $E_n$ are called the $n$-th Euler numbers.

From (6) and Theorem 2, we note that
\begin{eqnarray}
\tilde E_{n,q} (x) &= &\int_{\Bbb Z_p} (x+y)^{n} d\mu_{-q}(y) \notag
\\&=&
\sum_{l=0}^{n} \binom {n}{l} \int_{\Bbb Z_p} y^{l} d\mu_{-q}(y) x^{n-l}
\\&=& \sum_{l=0}^{n}  \binom {n}{l} \tilde E_{n,q} x^{n-l} \notag
\\&=&  {\left ( x+ \tilde E_q \right)}^n ,
 \quad \notag
\end{eqnarray}
where the usual convention about replacing $({\tilde E_q})^{l}$ by ${\tilde E}_{l,q}$.
By Theorem 2 and Theorem 4, we get

\begin{eqnarray}
 q \tilde E_{n,q} (1) + \tilde E_{n,q}  = \left \{
\begin{array}{cl}
 [2]_q , \  \ & \hbox{if} \ \ n=0,
\vspace{2mm} \\ 0, \ \ & \hbox{if} \ \ n>0.
\end{array}\right.
\end{eqnarray}

From (14) and (15), we have

\begin{eqnarray}
 q {\left( \tilde E_{q} +1 \right)}^n + \tilde E_{n,q}  = \left \{
\begin{array}{cl}
 [2]_q , \  \ & \hbox{if} \ \ n=0,
\vspace{2mm} \\ 0, \ \ & \hbox{if} \ \ n>0.
\end{array}\right.
\end{eqnarray}

For $n \in \Bbb N$, by (14) and (16), we have
\begin{eqnarray}
q^{2} \tilde E_{n,q} (2)
&=& q^{2} {\left( \tilde E_{q} +1 +1 \right)}^n  \notag \\
&=& q^2 \sum_{l=1}^{n} \binom {n}{l} {\left( \tilde E_{q} +1 \right)}^l + q {\left( 1+q - \tilde E_{0,q} \right)} \\
&=& q+ q^2 -q \sum_{l=0}^{n} \binom {n}{l} \tilde E_{l,q}  \notag \\
&=&  q+ q^2 -q  {\left (  \tilde E_{q} +1 \right )}^n \notag \\
&=&  q+q^2 + \tilde E_{n,q}. \notag
\end{eqnarray}

Therefore, by (17), we obtain the following theorem.
\begin{theorem}
For $n \in \Bbb N$, we have
$$
q^{2} \tilde E_{n,q} (2) = q + q^{2} + \tilde E_{n,q}.
$$
\end{theorem}

For $n \in \Bbb Z_+$, we have
\begin{eqnarray}
\tilde E_{n,q^{-1}} (1-x) &= &  \int_{\Bbb Z_p} (1-x+x_1 )^{n} d\mu_{-q^{-1}}(x_1 ) \notag
\\&=& (-1)^n \int_{\Bbb Z_p} (x_{1} + x)^{n} d\mu_{-q}(x_1)
\\&=&  (-1)^n \tilde E_{n,q} (x).  \notag
\end{eqnarray}

Therefore, by (18), we obtain the following theorem.
\begin{theorem}
For $n \in \Bbb Z_+$, we have
$$
\tilde E_{n,q^{-1}} (1-x) = (-1)^n  \tilde E_{n,q}(x).
$$
\end{theorem}

From (14), we have

\begin{eqnarray}
\int_{\Bbb Z_p} (1-x )^{n} d\mu_{-q}(x) \notag
&=& (-1)^n \int_{\Bbb Z_p} (x-1)^{n} d\mu_{-q}(x)
\\&=&  (-1)^n \tilde E_{n,q} (-1).
\end{eqnarray}

By Therorem 6 and (19), we get

\begin{eqnarray}
\int_{\Bbb Z_p} (1-x )^{n} d\mu_{-q}(x)
&=& \tilde E_{n,q^{-1}} (2) = 1 + q + q^2 \tilde E_{n,q^{-1}} \quad \text {if} \quad  n>0.
\end{eqnarray}

Therefore, by (20), we obtain the following theorem.
\begin{theorem}
For $n \in \Bbb N$, we have
$$
\int_{\Bbb Z_p} (1-x )^{n} d\mu_{-q}(x)
=  1 + q +q^2  \tilde E_{n,q^{-1}}.
$$
\end{theorem}

Let $C(\Bbb Z_p )$ be the space of continuous functions on $\Bbb Z_p$. For $f \in C(\Bbb Z_p )$, $p$-adic
analogue of Bernstein operator of order $n$ for $f$ is given by
\begin{eqnarray}
\Bbb B_n (f|x) &=& \sum_{k=0}^{n} B_{k,n} (x) f {\left( \frac {k}{n} \right)} \\
&=&  \sum_{k=0}^{n} f {\left( \frac {k}{n} \right)} \binom {n}{k} x^k (1-x)^{n-k}, \notag
\end{eqnarray}
where $n,k \in \Bbb Z_+$ (see [1,6,7]).

For $n,k \in \Bbb Z_+$, $p$-adic
Bernstein polynomials of degree $n$ is defined by
\begin{eqnarray}
B_{k,n }(x) &=& \binom {n}{k} x^k (1-x)^{n-k}, \quad x \in \Bbb Z_p \quad \text{(see [1,6,7]).}
\end{eqnarray}

Let us take the fermionic $p$-adic $q$-integral on $\Bbb Z_p$ for one Bernstein polynomials in (22) as
follows:

\begin{eqnarray}
\int_{\Bbb Z_p} { B_{k,n} (x )} d\mu_{-q}(x)&=& \binom {n}{k} \int_{\Bbb Z_p} x^{k}(1-x)^{n-k} d\mu_{-q}(x)\notag  \\
&=&  \binom {n}{k} \sum_{l=0}^{n-k} \binom {n-k}{l} (-1)^l  \int_{\Bbb Z_p} x^{k+l}d\mu_{-q}(x) \\
&=&  \binom {n}{k} \sum_{l=0}^{n-k} \binom {n-k}{l} (-1)^l  \tilde E_{k+l,q}  \notag
\end{eqnarray}

By simple calculation, we easily get

\begin{eqnarray}
\int_{\Bbb Z_p} { B_{k,n} (x )} d\mu_{-q}(x)&=& \int_{\Bbb Z_p} { B_{n-k,n} (1-x )} d\mu_{-q}(x) \notag \\
&=& \binom {n}{k} \sum_{l=0}^{k} \binom {k}{l} (-1)^{k+l} \int_{\Bbb Z_p} (1-x)^{n-l} d\mu_{-q}(x) \notag \\
&=&  \binom {n}{k} \sum_{l=0}^{k} \binom {k}{l} (-1)^{k+l} {\left ( 1 + q + q^2 \tilde E_{n-l, q^{-1}} \right)}  \\
&=&  \binom {n}{k} \sum_{l=0}^{k} \binom {k}{l} (-1)^{k+l}  q^2 \tilde E_{n-l,q^{-1}} \quad \text {if} \quad  n > k. \notag
\end{eqnarray}

Therefore, by (23) and (24), we obtain the following theorem.

\begin{theorem}
For $n \in \Bbb Z_+$ with $n > k$, we have
$$
\sum_{l=0}^{n-k} \binom {n-k}{l} (-1)^l \tilde E_{k+l,q}=
\sum_{l=0}^{k} \binom {k}{l} (-1)^{k+l} q^2 \tilde E_{n-l,q^{-1}}.
$$
\end{theorem}
In particular, $k=0$, we get
$$
\sum_{l=0}^{n} \binom {n}{l} (-1)^l \tilde E_{l,q}= q^2
\tilde E_{n,q^{-1}}.
$$

By Theorem 1 and Theorem 2, we get
$$
\sum_{l=0}^{n-k} \binom {n-k}{l} (-1)^l H_{k+l} (-q^{-1})=
\sum_{l=0}^{k} \binom {k}{l} (-1)^{k+l} q^2 H_{n-l} (-q),
$$
where $H_n (-q)$ are called the $n$-th Frobenius-Euler numbers.

\end{document}